\documentclass[11pt]{article}
\usepackage[T1]{fontenc}
\usepackage{lmodern}
\usepackage[margin=1in]{geometry}
\usepackage{amsmath,amssymb,amsthm,mathtools}
\usepackage{microtype}
\usepackage{indentfirst}
\usepackage{enumitem}
\usepackage{needspace}
\setlist[itemize]{leftmargin=*,itemsep=3pt,topsep=5pt}
\usepackage{tikz}
\usepackage{tikz-cd}
\usetikzlibrary{arrows.meta}
\usepackage{caption}
\usepackage{booktabs}
\usepackage{placeins}
\usepackage[colorlinks=true,linkcolor=blue!45!black,citecolor=blue!45!black,urlcolor=blue!45!black,bookmarksnumbered=true]{hyperref}
\usepackage{etoolbox}
\makeatletter
\patchcmd{\l@section}{\bfseries}{\normalfont}{}{}
\makeatother
\numberwithin{equation}{section}
\newtheorem{theorem}{Theorem}[section]
\newtheorem{proposition}[theorem]{Proposition}
\newtheorem{lemma}[theorem]{Lemma}
\newtheorem{corollary}[theorem]{Corollary}
\theoremstyle{definition}
\newtheorem{definition}[theorem]{Definition}
\newtheorem{example}[theorem]{Example}
\theoremstyle{remark}

\newcommand{\C}{\mathbb C}
\newcommand{\R}{\mathbb R}
\newcommand{\Z}{\mathbb Z}
\newcommand{\T}{\mathbb T}
\newcommand{\Shape}{\mathcal S}
\newcommand{\Vol}{\operatorname{Vol}}
\newcommand{\Area}{\operatorname{Area}}
\newcommand{\conv}{\operatorname{conv}}

\newcommand{\id}{\mathrm{id}}
\DeclareMathOperator{\Int}{Int}

\newcommand{\address}[1]{\AtEndDocument{\par\bigskip{\small\noindent #1\par}}}
\address{Beijing Institute of Mathematical Sciences and Applications (BIMSA), Beijing, 101408, China; \\ E-mail: \textit{yuanhang@bimsa.cn}}

\title{An exact K\"ahler form on $(\C^*)^2$ without a global SYZ fibration}
\author{Hang Yuan}

\date{}

\begin{document}
\maketitle
\begin{abstract}
We construct an explicit exact K\"ahler form on \((\mathbb C^*)^2\) admitting no almost-toric fibration and hence no smooth Lagrangian torus fibration. For the logarithmic holomorphic volume form, we also rule out any special Lagrangian fibrations with compact fibers and dense regular locus, even when arbitrary compact connected singular fibers are allowed. For higher dimensions, taking products also yields explicit exact K\"ahler forms on \((\mathbb C^*)^n\), for every \(n\ge2\), admitting no smooth proper Lagrangian torus fibration.
\end{abstract}

\section{Introduction}\label{sec:introduction}

The Strominger--Yau--Zaslow (SYZ) picture relates mirror symmetry to dual
special Lagrangian fibrations, whose regular fibers are tori
\cite{SYZ1996}. Over the regular locus, the dual torus parametrizes
rank-one unitary local systems on the
original fiber. In the log Calabi--Yau setting, this leads to studying
such fibrations on anticanonical complements equipped with their
holomorphic volume forms \cite[Sections~1--2]{Auroux2007}.

A basic model is the algebraic torus $X=(\C^*)^2$, the complement of the
toric boundary $D$ in $\overline X=\mathbb P^1\times\mathbb P^1$.
Since $K_{\overline X}+D\sim0$, $X$ is a log Calabi--Yau surface.
Its standard logarithmic volume form is
\begin{equation}\label{eq:log-volume-form}
\Omega=\frac{dz_1}{z_1}\wedge\frac{dz_2}{z_2}.
\end{equation}
Write $\T^n=(\R/\Z)^n$, with each angular coordinate of period one.
The logarithmic polar coordinates $z_j=e^{a_j+2\pi i\theta_j}$ identify
$X$ diffeomorphically with $T^*\T^2$. For the
standard cotangent K\"ahler form, the logarithm map
$\operatorname{Log}(z_1,z_2)=(\log|z_1|,\log|z_2|)$ is a proper
submersion whose torus fibers are special Lagrangian with respect to
$\Omega$.

We ask whether a special Lagrangian fibration with compact fibers must
still exist when the K\"ahler form changes. The logarithm fibers may cease
to be Lagrangian, but other fibers could replace them. In this paper, we give an
explicit K\"ahler form for which no such fibration exists.

With the conventions
$d^c=\frac{i}{4\pi}(\bar\partial-\partial)$ and
$dd^c=\frac{i}{2\pi}\partial\bar\partial$, set
\begin{equation}\label{eq:potential}
\Phi(z_1,z_2)=\log\!\left(1+|z_1|^2+|z_2|^2+\tfrac14|z_1^2+z_2^2|^2\right),
\quad \lambda=d^c\Phi,\quad \omega=d\lambda=dd^c\Phi.
\end{equation}

The form $\omega$ is a real-analytic exact K\"ahler form on $X$ with
$\int_X\omega^2/2=1$, as shown in Proposition~\ref{prop:kahler-volume}.
For the fixed logarithmic volume form $\Omega$, we classify all compact
special Lagrangian surfaces.
For the diagonal circle action
$\tau\cdot(z_1,z_2)=(e^{2\pi i\tau}z_1,e^{2\pi i\tau}z_2)$,
$\tau\in\R/\Z$, we normalize the moment map $\mu$ by
\begin{equation}\label{eq:moment}
\mu=
\frac{|z_1|^2+|z_2|^2+\tfrac12|z_1^2+z_2^2|^2}
{1+|z_1|^2+|z_2|^2+\tfrac14|z_1^2+z_2^2|^2}.
\end{equation}

Define the joint-level map
\begin{equation}\label{eq:joint-level-map}
F=(\mu,\log|z_2/z_1|):X\longrightarrow\R^2.
\end{equation}

\Needspace{9\baselineskip}
\begin{theorem}[Classification of compact special Lagrangians]\label{thm:compact-special-lagrangian}
The compact connected embedded special Lagrangian surfaces without boundary
in $(X,\omega,\Omega)$ are exactly the tori
\[
L_{c,R}=\{\mu=c,\ |z_2/z_1|=R\},
\]
where $R>0$, with $0<c<2$ if $R\ne1$ and $0<c<1$ if $R=1$.
These tori are the fibers of the smooth proper special Lagrangian torus
fibration $F|_U:U\to B_0$, where $U=F^{-1}(B_0)$ and explicitly
\[
U=X\setminus\bigl\{(z_1,z_2)\in X:
|z_1|=|z_2|,\ |z_1^2+z_2^2|\geq2\bigr\},\quad
B_0=\bigl((0,2)\times\R\bigr)\setminus\bigl([1,2)\times\{0\}\bigr).
\]
The set $U$ is open and dense in $X$, and the base $B_0$ is diffeomorphic
to $\R^2$.
\end{theorem}

The tori in Theorem~\ref{thm:compact-special-lagrangian} already fiber an open
dense subset of $X$. To extend this picture to all of $X$, we must allow
singular fibers, and a natural requirement is that these fibers be compact
and connected. We make these requirements precise in the following definition.

\begin{definition}\label{def:special-fibration}
A \emph{special Lagrangian fibration} of $(X,\omega,\Omega)$ is a
continuous surjection $\pi:X\to B$ with compact connected fibers, where
$B$ is a connected Hausdorff second-countable surface. For some open subset $B_{\mathrm{reg}}\subset\Int(B)$, $X_{\mathrm{reg}}:=\pi^{-1}(B_{\mathrm{reg}})$
is dense in $X$ and $\pi|_{X_{\mathrm{reg}}}$ is a smooth submersion with
special Lagrangian torus fibers.
No further conditions are imposed on the fibers over $B\setminus B_{\mathrm{reg}}$.
\end{definition}

The above classification result and a continuity argument using the noncompact joint
levels of $F$ give the following global obstruction.

\Needspace{4\baselineskip}
\begin{corollary}[Special Lagrangian obstruction]\label{cor:special-lagrangian}
The triple $(X,\omega,\Omega)$ admits no special Lagrangian
fibration with compact fibers.
No further restriction is imposed on the singular fibers.
\end{corollary}

If compactness of fibers is dropped, a naive completion is available.
Indeed, fix $0<a<1$ and identify the domain
$B_0\setminus([a,1)\times\{0\})$
diffeomorphically with $S^2\setminus\{N\}$. Composing $F|_U$ with this
identification on the corresponding tori, and sending all remaining
points of $X$ to the north pole $N$, gives a continuous surjection $\pi:X\to S^2$
with a dense regular locus of special Lagrangian torus fibers. Its sole
exceptional fiber is
\[
\pi^{-1}(N)=\{|z_2/z_1|=1,\ \mu\ge a\}
\cong\T^2\times[0,\infty),
\]
which is connected but noncompact. Thus this map satisfies every
condition in Definition~\ref{def:special-fibration} except compactness of the
exceptional fiber.

The preceding obstruction uses the fixed holomorphic volume form \(\Omega\). What restrictions does \(\omega\) alone impose on Lagrangian torus fibrations?
We next give a purely symplectic obstruction:

\Needspace{4\baselineskip}
\begin{theorem}\label{thm:main}
The symplectic manifold $(X,\omega)$ admits no almost-toric Lagrangian
fibration and hence no smooth Lagrangian torus fibration.
\end{theorem}

The regularity criterion in Section~\ref{sec:dimension-four} shows that any
almost-toric fibration on $(X,\omega)$ must be a smooth Lagrangian torus
fibration. Then, we compare symplectic volume with Sikorav--Eliashberg's
Lagrangian shape \cite{Sikorav1989,Sikorav1991,Eliashberg}:
using Gromov's exact Lagrangian intersection theorem \cite{Gromov}
and Duistermaat's global action--angle coordinates \cite{Duistermaat},
Section~\ref{sec:shape} establishes a shape--volume identity under suitable
hypotheses. In Section~\ref{sec:tori}, we construct the required exact symplectic
embedding and use contraction centers and the action classes of the
special Lagrangian tori to show that the shape has area at least $4/3$,
contradicting the identity since the symplectic volume is $1$.
Taking products with normalized Fubini--Study factors, the shape--volume
argument will also yield the following:

\begin{corollary}[Higher-dimensional cases]\label{cor:higher-dimensional}
For every $n\ge2$, the algebraic torus $(\C^*)^n$ admits a real-analytic
exact K\"ahler form admitting no smooth proper Lagrangian torus fibration.
\end{corollary}

The two obstructions have complementary scopes.
Corollary~\ref{cor:special-lagrangian} concerns the fixed logarithmic volume form
and allows arbitrary compact connected singular fibers.
Theorem~\ref{thm:main} is independent of the holomorphic volume form but
restricts the singularities to the almost-toric class. General Lagrangian
fibrations with singularities outside this class are not addressed here.

\paragraph*{\normalfont\itshape Significance and scope.}

The underlying algebraic log Calabi--Yau structure,
K\"ahler class, and total volume do not determine the existence of a special Lagrangian fibration.
Indeed, the product Fubini--Study form $\omega_{\mathrm{FS}}$ satisfies
$[\omega_{\mathrm{FS}}]=[\omega]=0$ and has the same total volume $1$,
yet admits the special Lagrangian logarithm fibration for the same
$\Omega$. 
Note that our metric is incomplete and positive K\"ahler--Einstein. In particular, it is not Ricci-flat.
Our example shows that the special Lagrangian fibration picture of SYZ
does not extend automatically to arbitrary K\"ahler metrics on log
Calabi--Yau manifolds. By contrast, existence results are available for
complete Tian--Yau metrics. Indeed, Collins--Jacob--Lin construct special Lagrangian
fibrations on complements of smooth anticanonical divisors in del Pezzo
and rational elliptic surfaces \cite{CollinsJacobLin2021}.

Our results suggest that the singularities of SYZ fibrations may be more
complicated than we had expected. The SYZ construction begins by dualizing
the smooth torus fibers, which typically cover an open dense subset, and
gluing the resulting mirror charts using Floer-theoretic
corrections~\cite{YuanThesis}.
One might hope to recover the missing parts of the mirror from a sufficiently developed Floer
theory of the singular fibers. But, our result suggests that there may be a more basic difficulty: the regular fibers may not admit a completion by compact singular fibers at all.

When constructing the SYZ mirror of an open log Calabi--Yau manifold, it may
be \textit{necessary} to remember how the manifold $X$ sits inside a compactification $Y$, as our results suggest.
In our example, \(\omega\) is not the restriction of any K\"ahler form on the standard toric compactification \(\mathbb P^1\times\mathbb P^1\).
More generally, let $(Y,\omega_Y)$
be a compact K\"ahler manifold with an anticanonical divisor $D$.
Suppose that $Y\setminus D$, equipped with the restriction of $\omega_Y$,
admits a Lagrangian torus fibration, and let $L$ be a regular fiber.
A holomorphic disk
$u:(\mathbb D,\partial\mathbb D)\to(Y,L)$ can meet $D$ in its
interior. In the family Floer SYZ construction of the Landau--Ginzburg
superpotential, the weight of a counted disk depends on its area
$\int_{\mathbb D}u^*\omega_Y$. The author's thesis~\cite{YuanThesis}
shows that the local disk potentials glue to a globally defined mirror superpotential $W$. For the $A_n$ smoothing, disk counts in a (partial)
compactification give local expressions for $W$, and matching these local expressions
on overlaps is a key step in determining its
SYZ mirror~\cite[Section~3.3]{YuanAnSYZ}. The toric Calabi-Yau,
conifold, and Grassmannian examples follow the same pattern of using
disk potentials to relate the boundary geometry to the family Floer SYZ
mirror~\cite{YuanDiskCounting,YuanLocalSYZ,YuanConifold,YuYuanGrassmannian}.

\paragraph*{\normalfont\itshape Acknowledgment and AI disclosure.}
The project was initially motivated in part by the role of compactification in the author’s earlier work on how disk counts and mirror superpotentials help determine SYZ mirrors. These considerations prompted a discussion with AI about how the choice of compactification might bear on the existence of a Lagrangian torus fibration on the complement, and that discussion ultimately led to the present paper. 
The author acknowledges substantial use of AI tools in preparing and editing the manuscript.

\section{Construction of the K\"ahler form}\label{sec:geometry}

Let $X=(\C^*)^2$ with $\Phi$, $\lambda=d^c\Phi$, and
$\omega=d\lambda$ as in \textup{(\ref{eq:potential})}.
Writing $z_j=e^{a_j+2\pi i\theta_j}$, $j=1,2$, defines logarithmic
coordinates $(a,\theta)\in\R^2\times\T^2$ and identifies $X$
diffeomorphically with $T^*\T^2$. In particular, $X\simeq\T^2$.
With our normalization of $d^c$, we have $d^c\log|z_j|^2=d\theta_j$.

Consider the unitary linear map
\begin{equation}\label{eq:unitary}
T=(u,v):\C^2\longrightarrow\C^2,
\qquad
T(z_1,z_2)=(u(z_1,z_2),v(z_1,z_2))
=\left(\frac{z_1+i z_2}{\sqrt2},\frac{z_1-i z_2}{\sqrt2}\right).
\end{equation}
Since $z_1=(u+v)/\sqrt2$ and $z_2=(u-v)/(i\sqrt2)$, its restriction is a biholomorphism
\[
T:(\C^*)^2\xrightarrow{\sim}
\C^2\setminus\bigl(\{u=v\}\cup\{u=-v\}\bigr).
\]
Since $|u|^2+|v|^2=|z_1|^2+|z_2|^2$ and $uv=(z_1^2+z_2^2)/2$, we have
\begin{equation}\label{eq:factorization}
\Phi(z_1,z_2)=\log(1+|u|^2)+\log(1+|v|^2).
\end{equation}

\begin{proposition}\label{prop:kahler-volume}
The form $\omega$ is a real-analytic exact K\"ahler form on $X$, and
$\int_X\omega^2/2=1$.
\end{proposition}

\begin{proof}
The expression in \textup{(\ref{eq:factorization})} defines a real-analytic potential on all of
$\C^2$, with positive definite complex Hessian
$\operatorname{diag} ((1+|u|^2)^{-2},(1+|v|^2)^{-2} )$
in the $(u,v)$ coordinates. Thus $\omega=dd^c\Phi$ is a real-analytic
K\"ahler form on $\C^2$, whose restriction to $X$ is exact with primitive
$\lambda=d^c\Phi$.

To compute its volume, put
\[
\sigma=dd^c\log(1+|w|^2)
=\frac{i\,dw\wedge d\bar w}{2\pi(1+|w|^2)^2}
\quad\text{on }\C.
\]
We have
\[
\int_{\C}\sigma
=\int_0^\infty\frac{2r}{(1+r^2)^2}\,dr=1.
\]
Let $\operatorname{pr}_j$ be the coordinate projections on the target $\C^2$ of $T$.
By \textup{(\ref{eq:factorization})},
$
\omega=T^*\!\left(\operatorname{pr}_1^*\sigma+\operatorname{pr}_2^*\sigma\right)$ and $
\omega^2 / 2
=T^*\! (\operatorname{pr}_1^*\sigma\wedge\operatorname{pr}_2^*\sigma )$.
The complement of $T(X)$ consists of the two complex lines $\{u=v\}$ and
$\{u=-v\}$, which have zero four-dimensional measure. Thus change of variables
by $T$ and Fubini's theorem give
\[
\int_X\frac{\omega^2}{2}
=\int_{T(X)}\operatorname{pr}_1^*\sigma\wedge\operatorname{pr}_2^*\sigma=\int_{\C^2}\operatorname{pr}_1^*\sigma\wedge\operatorname{pr}_2^*\sigma
=\left(\int_{\C}\sigma\right)^2=1. \qedhere
\]
\end{proof}

\section{Special Lagrangian tori and fibrations}\label{sec:special-lagrangian}

An \emph{almost Calabi--Yau surface} $(M,\omega_M,\Omega_M)$ is a
K\"ahler surface equipped with a nowhere-vanishing holomorphic two-form
$\Omega_M$; the metric need not be Ricci-flat.
For an oriented embedded Lagrangian surface $L\subset M$, define its
\emph{phase map} by
\[
\vartheta_L=\arg(\Omega_M|_L):L\longrightarrow\R/2\pi\Z.
\]
Here $\vartheta_L(x)=\arg\bigl(\Omega_M(v_1,v_2)\bigr)$ for any
positively oriented basis $(v_1,v_2)$ of $T_xL$; this is independent of
the basis since a change of such basis multiplies the value by a positive
real number.
The surface is \emph{special Lagrangian of phase $\vartheta$} if
$\vartheta_L\equiv\vartheta$.
Equivalently, $(e^{-i\vartheta}\Omega_M)|_L$ is a positive real area form.
In particular, $
\operatorname{Im}(e^{-i\vartheta}\Omega_M)|_L=0$.

We first record how Hamiltonian symmetries constrain compact special
Lagrangian surfaces.

\begin{lemma}[Hamiltonian symmetry]\label{lem:sl-hamiltonian-symmetry}
Let $(M,\omega_M,\Omega_M)$ be an almost Calabi--Yau surface, and let $Y$
be a vector field satisfying
\[
\iota_Y\omega_M=-dH,\qquad \mathcal L_Y\Omega_M=0
\]
for a smooth real-valued function $H$.
If $L\subset M$ is a closed connected special Lagrangian surface, then
$H|_L$ is constant and $Y$ is tangent to $L$.
\end{lemma}

\begin{proof}
Let $J$ be the complex structure of $M$, and let $dA_g$ be the oriented
area form on $L$ induced by $g=\omega_M(\cdot,J\cdot)$.
Set $\Psi=e^{-i\vartheta_L}\Omega_M$, $ h=H|_L$, and $
\beta=(\iota_Y\operatorname{Im}\Psi)|_L$.
The special Lagrangian condition gives $\Psi|_L=\rho\,dA_g$ for a
smooth function $\rho:L\to(0,\infty)$.

Along $L$, decompose $Y|_L=Y^\top+Y^\perp$ into its $g$-orthogonal
tangent and normal components. For $u,v\in TL$, the Lagrangian condition
gives $g(Ju,v)=\omega_M(u,v)=0$. Thus $J(TL)$ is the normal bundle. Define $W:=-JY^\perp$; this is a tangent
vector field on $L$ and $JW=Y^\perp$ as $J^2=-\id$.
For $v\in TL$, we have 
$dh(v)=-\omega_M(Y,v)=-\omega_M(JW,v)=g(W,v)$,
so $W=\nabla_Lh$ for the induced metric on $L$.
For $v\in TL$, both $Y^\top$ and $v$ are tangent to $L$.
Since $\Psi|_L=\rho\,dA_g$ is real, we have
$\operatorname{Im}\Psi(Y^\top,v)=0$. Now, we compute
\[
\begin{aligned}
\beta(v)&=\operatorname{Im}\Psi(Y,v)
=\operatorname{Im}\Psi(Y^\perp,v)
=\operatorname{Im}\Psi(JW,v) \\
&=\operatorname{Im}\bigl(i\Psi(W,v)\bigr)
=\operatorname{Re}\Psi(W,v)
=\rho\,dA_g(\nabla_Lh,v).
\end{aligned}
\]

The holomorphic top-degree form $\Psi$ is closed, and
$\mathcal L_Y\Psi=0$ by hypothesis. Cartan's formula gives $d\beta=0$,
so Stokes' theorem on the closed surface $L$ yields
\[
0=\int_Ld(h\beta)=\int_Ldh\wedge\beta
=\int_L\rho\,|\nabla_Lh|_g^2\,dA_g.
\]
Positivity of $\rho$ implies $\nabla_Lh=0$. Thus $h=H|_L$ is constant
by connectedness, and $Y^\perp=J\nabla_Lh=0$ shows that $Y$ is
tangent to $L$.
\end{proof}

We now apply this lemma to $(X,\omega,\Omega)$.
The diagonal circle action
$\tau\cdot(z_1,z_2)=(e^{2\pi i\tau}z_1,e^{2\pi i\tau}z_2)$, $\tau\in\R/\Z$,
preserves $\Phi$ and $\lambda$. If $Y$ is its period-one infinitesimal generator, then
the function $\mu$ from \textup{(\ref{eq:moment})} satisfies $\mu=\lambda(Y)$.
Invariance of $\lambda$ implies $\iota_Y\omega=-d\mu$.
Set $R:X\to (0,\infty)$, $R(z_1,z_2)=|z_2/z_1|$.
The circle action also preserves $\Omega$, and
\begin{equation}\label{eq:sl-circle-contraction}
\iota_Y\omega=-d\mu,\qquad
\iota_Y\Omega=2\pi i\,d\log(z_2/z_1),\qquad
\iota_Y\operatorname{Im}\Omega=2\pi\,d(\log R).
\end{equation}

\begin{lemma}\label{lem:sl-constant-levels}
Every closed connected special Lagrangian surface $L$ in $(X,\omega,\Omega)$
satisfies $\vartheta_L\equiv0$ or $\vartheta_L\equiv\pi$, and both
$\mu$ and $R$ are constant on $L$.
\end{lemma}

\begin{proof}
The group $H_2(X;\Z)\cong\Z$ is generated by the standard product torus $T$,
on which the period $\int_T\Omega$ is $-4\pi^2$.
Thus $\int_L\Omega\in-4\pi^2\Z$. Since $\vartheta_L$ is constant
and $(e^{-i\vartheta_L}\Omega)|_L$ is a positive area form, we must have $\vartheta_L\equiv0$ or $\vartheta_L\equiv\pi$.
By Lemma~\ref{lem:sl-hamiltonian-symmetry} and \textup{(\ref{eq:sl-circle-contraction})}, $\mu$ is constant and $Y$ is
tangent to $L$. Since $\operatorname{Im}\Omega|_L=0$, the last identity
in \textup{(\ref{eq:sl-circle-contraction})} gives $d(\log R)|_L=0$.
Thus $R$ is constant on $L$.
\end{proof}

\begin{proposition}[Geometry of the joint levels]\label{prop:explicit-tori}
The joint-level map $F$ defined in \textup{(\ref{eq:joint-level-map})} is a submersion with image
$F(X)=(0,2)\times\R$.
For $c\in\R$ and $R>0$, its joint levels
\[
L_{c,R}=F^{-1}(c,\log R)=\{\mu=c,\ |z_2/z_1|=R\}
\]
have the following diffeomorphism types:
\[
L_{c,R}\cong
\begin{cases}
\T^2, & R\ne1,\quad 0<c<2,\\
\T^2, & R=1,\quad 0<c<1,\\
(S^1\times\R)\sqcup(S^1\times\R), & R=1,\quad 1\leq c<2,\\
\varnothing, & c\leq0\ \text{or}\ c\geq2.
\end{cases}
\]
Every nonempty joint level is an embedded special Lagrangian surface,
admitting an orientation of phase $\pi$. For every torus level, the inclusion in $X$ induces an isomorphism
on integral first homology.
\end{proposition}

\begin{proof}
Write $z_j=|z_j|e^{2\pi i\theta_j}$, and set
$r=|z_1|^2$ and $\varphi=\theta_2-\theta_1\in\R/\Z$.
For $R>0$, put
$k=1+R^2$ and $h(\varphi)=\frac14\left|1+R^2e^{4\pi i\varphi}\right|^2$.
Since $z_2/z_1=Re^{2\pi i\varphi}$, equation \textup{(\ref{eq:moment})} becomes
\begin{equation}\label{eq:radial-moment}
\mu=\frac{kr+2hr^2}{1+kr+hr^2}.
\end{equation}
Here $k>0$ and $h\geq0$, so $0<\mu<2$ and
\[
\frac{\partial\mu}{\partial r}
=
\frac{k+4hr+khr^2}{(1+kr+hr^2)^2}>0.
\]
The Jacobian minor of $F=(\mu,\log R)$ in the $(r,R)$ variables has
determinant $R^{-1}\partial\mu/\partial r>0$. Thus $F$ is a submersion,
every joint level is a smooth embedded surface, and
$L_{c,R}=\varnothing$ for $c\notin(0,2)$.

Fix $0<c<2$ and $R>0$, and let $A_{c,R}\subset\R/\Z$ be the set of
relative angles at which $\mu(r,R,\varphi)=c$ has a solution $r>0$.
If \(R\ne1\), then \(h>0\). For each fixed \(\varphi\), the right-hand side of \textup{(\ref{eq:radial-moment})} increases from \(0\) to \(2\) as \(r\) increases from \(0\) to \(+\infty\).
If $R=1$, then $h(\varphi)=\cos^2(2\pi\varphi)$.
Away from its two zeros $\varphi=1/4,3/4$, the radial expression of $\mu$ still
has range $(0,2)$; at the two zeros, it is $2r/(1+2r)$, with range $(0,1)$.
Consequently,
\[
A_{c,R}=
\begin{cases}
\R/\Z, & R\ne1\ \text{or}\ 0<c<1,\\
(\R/\Z)\setminus\{1/4,3/4\}, & R=1,\quad 1\leq c<2.
\end{cases}
\]
In particular, $A_{c,R}$ is nonempty for every $0<c<2$ and $R>0$.
Since $\log R$ ranges over $\R$, this proves $F(X)=(0,2)\times\R$.

For each \(\varphi\in A_{c,R}\), a positive solution of \(\mu(r,R,\varphi)=c\) exists by definition of \(A_{c,R}\), and it is unique because \(\partial\mu/\partial r>0\). Denote it by \(r_{c,R}(\varphi)\). The implicit function theorem shows that this solution depends smoothly on the parameters \((c,R,\varphi)\).
Define $f_{c,R}\colon(\R/\Z)\times A_{c,R}\to X$ by
\begin{equation}\label{eq:torus-parameterization}
f_{c,R}(\theta_1,\varphi)
=\left(\sqrt{r_{c,R}(\varphi)}e^{2\pi i\theta_1} \ , \
R\sqrt{r_{c,R}(\varphi)}e^{2\pi i(\theta_1+\varphi)}\right).
\end{equation}
One can check that $f_{c,R}$ is a
diffeomorphism onto $L_{c,R}$, and
$L_{c,R}\cong(\R/\Z)\times A_{c,R}$.
When $A_{c,R}=\R/\Z$, this gives $L_{c,R}\cong\T^2$.
Otherwise, $A_{c,R}$ is the disjoint union of two open intervals, so
$L_{c,R}$ consists of two components, each diffeomorphic to $S^1\times\R$.

To check the Lagrangian condition on any nonempty level, note that varying
$\theta_1$ while holding $\varphi$ fixed gives the diagonal circle action.
Thus $df_{c,R}(\partial_{\theta_1})=Y$. Set
$V=df_{c,R}(\partial_\varphi)$. Since $\mu\circ f_{c,R}=c$, the identity
$\iota_Y\omega=-d\mu$ gives
$\omega(Y,V)=-d\mu(V)=-\partial_\varphi(\mu\circ f_{c,R})=0$.
The vectors $Y$ and $V$ span the tangent plane of $L_{c,R}$ at every point,
so this two-dimensional surface is Lagrangian.

For the phase, write $r=r_{c,R}(\varphi)$ and $r'=dr/d\varphi$.
Since $R$ is constant, a direct computation gives
\[
\alpha:=f_{c,R}^*\frac{dz_1}{z_1}
=\frac{r'}{2r}\,d\varphi+2\pi i\,d\theta_1,
\qquad
f_{c,R}^*\frac{dz_2}{z_2}=\alpha+2\pi i\,d\varphi.
\]
We then obtain
\[
f_{c,R}^*\Omega
=\alpha\wedge(\alpha+2\pi i\,d\varphi)
=2\pi i\left(\frac{r'}{2r}\,d\varphi+2\pi i\,d\theta_1\right)\wedge d\varphi =-4\pi^2\,d\theta_1\wedge d\varphi.
\]
With the orientation induced by $d\theta_1\wedge d\varphi$, this gives
$\vartheta_{L_{c,R}}\equiv\pi$. Hence every nonempty level is also special
Lagrangian.
Finally, for a torus level, $A_{c,R}=\R/\Z$. The coordinate cycles in the
$(\theta_1,\varphi)$ parametrization have angular winding vectors $(1,1)$
and $(0,1)$, respectively. Their winding matrix has determinant one, so
the inclusion induces an isomorphism on integral first homology.
\end{proof}

In the following, we continue to use the coordinates and notation from the proof of Proposition~\ref{prop:explicit-tori}.

\begin{proof}[Proof of Theorem~\ref{thm:compact-special-lagrangian}]
Let $L\subset X$ be a closed connected embedded special Lagrangian surface.
By Lemma~\ref{lem:sl-constant-levels}, there are constants $c\in(0,2)$ and $R>0$
such that
$L\subset L_{c,R}$.
By Proposition~\ref{prop:explicit-tori}, $L_{c,R}$ is a smooth surface.
Since $L$ is connected, it is a connected component of $L_{c,R}$.
The classification in Proposition~\ref{prop:explicit-tori} shows that a compact
component exists exactly in the torus cases stated in Theorem~\ref{thm:compact-special-lagrangian};
in those cases the level is connected, so $L=L_{c,R}$.
Conversely, Proposition~\ref{prop:explicit-tori} shows that every torus in the list
is special Lagrangian. This completes the classification.

It remains to identify the union of these tori and verify properness.
By \textup{(\ref{eq:moment})},
\[
\mu-1=
\frac{\tfrac14|z_1^2+z_2^2|^2-1}
{1+|z_1|^2+|z_2|^2+\tfrac14|z_1^2+z_2^2|^2}.
\]
Consequently, $U=F^{-1}(B_0)=\{R\ne1\}\cup\{\mu<1\}$, which is open and dense
in $X$. The coordinate change
$(c,t)\mapsto(\log(c/(2-c)),t)$ identifies $B_0$ with the slit plane
$\R^2\setminus([0,\infty)\times\{0\})$. Polar coordinates identify
the latter with $(0,\infty)\times(0,2\pi)$, which is diffeomorphic to
$\R^2$. The preceding classification shows that
$F|_U:U\to B_0$ is a smooth surjective submersion whose
fibers are precisely all closed connected special Lagrangian surfaces in $X$.

Recall that the radial solution $r_{c,R}(\varphi)$ constructed above
depends smoothly on $(c,\log R,\varphi)\in B_0\times\R/\Z$
by the implicit function theorem. Thus the parametrizations
\textup{(\ref{eq:torus-parameterization})} give a diffeomorphism
\[
B_0\times\T^2\longrightarrow U,\qquad
\bigl((c,t),(\theta_1,\varphi)\bigr)
\longmapsto f_{c,e^t}(\theta_1,\varphi),
\]
under which $F|_U$ is projection to $B_0$. Since $\T^2$ is compact,
this projection is proper.
\end{proof}

The classification identifies all possible regular fibers of a special
Lagrangian fibration. We show that continuity forces any noncompact
joint level of $F$ to lie in a single fiber of the proposed fibration,
contradicting compactness.

\begin{proof}[Proof of Corollary~\ref{cor:special-lagrangian}]
Suppose that $\pi:X\to B$ is a special Lagrangian fibration in
the sense of Definition~\ref{def:special-fibration}. Choose
$b\in[1,2)\times\{0\}$. By Proposition~\ref{prop:explicit-tori}, the joint level
$F^{-1}(b)$ is nonempty and noncompact; it is closed in $X$ by continuity
of $F$. Fix any $p\in F^{-1}(b)$.
Since $X_{\mathrm{reg}}$ is dense, choose $p_k\in X_{\mathrm{reg}}$
with $p_k\to p$, and let $L_k=\pi^{-1}(\pi(p_k))$.
By Theorem~\ref{thm:compact-special-lagrangian}, each $L_k$ is the entire
joint level $F^{-1}(F(p_k))$.

Fix any $q\in F^{-1}(b)$.
By Proposition~\ref{prop:explicit-tori}, $F$ is a submersion, so there is a smooth
local section $s:V\to X$ on a neighborhood $V$ of $b$ such that
$
F\circ s=\id_V$ and $s(b)=q$.
Since $F(p_k)=b_k\to b$, for all sufficiently large $k$ we may set
$q_k=s(F(p_k))$. Then $q_k\to q$ and
$F(q_k)=F(p_k)$, so $q_k\in L_k$ and $\pi(q_k)=\pi(p_k)$.
Continuity of $\pi$ and the Hausdorff property of $B$ give
$\pi(q)=\pi(p)$. As $q$ was arbitrary, we obtain
\[
F^{-1}(b)\subseteq\pi^{-1}(\pi(p)).
\]
The left-hand side is closed in $X$ and noncompact, whereas the
right-hand side is compact, a contradiction.
\end{proof}


\section{Almost-toric fibrations}\label{sec:dimension-four}

We now turn to the symplectic obstruction in Theorem~\ref{thm:main}.
For exact symplectic four-manifolds homotopy equivalent to $\T^2$, we
show that every almost-toric fibration is smooth. This reduces the
problem to ruling out smooth Lagrangian torus fibrations, which we do
using Lagrangian shape in Sections~\ref{sec:shape} and~\ref{sec:tori}.

\begin{definition}[Smooth Lagrangian torus fibration]\label{def:smooth-fibration}
Let $(M^{2n},\omega_M)$ be a connected symplectic manifold without boundary.
A \emph{smooth Lagrangian torus fibration} is a smooth proper surjective
submersion $\pi:M\to B$ onto a connected smooth $n$-manifold without boundary,
with connected Lagrangian fibers. By the Arnold--Liouville theorem, every
fiber is diffeomorphic to $\T^n$.
\end{definition}

In dimension four, we also allow almost-toric singularities. The following
definition uses the local models of Leung and Symington
\cite[Definitions~2.1--2.2]{LeungSymington}.

\begin{definition}[Almost-toric fibration]\label{def:almost-toric-fibration}
Let $(M^4,\omega_M)$ be a connected symplectic manifold without boundary.
An \emph{almost-toric fibration} is a proper surjective map
$\pi:M\to B$ with connected fibers onto a connected surface $B$, possibly
with boundary and corners. In Darboux
coordinates $(x_1,y_1,x_2,y_2)$ with
$\omega_M=dx_1\wedge dy_1+dx_2\wedge dy_2$, and suitable coordinates on the
base, the map has one of the forms
\[
\begin{aligned}
(x_1,x_2) & &&\text{(regular)},\\
(x_1,x_2^2+y_2^2) & &&\text{(rank-one elliptic)},\\
(x_1^2+y_1^2,x_2^2+y_2^2) & &&\text{(elliptic--elliptic)},\\
(x_1y_1+x_2y_2,x_1y_2-x_2y_1) & &&\text{(focus--focus)}.
\end{aligned}
\]
Regular fibers are Lagrangian two-tori; fibers over boundary edges and
corners are circles and points, respectively. Interior singular fibers are
nodal and may contain several focus--focus points.
\end{definition}

Let $u,v$ be the component functions listed above. For each local model, write
\[
\pi:(x_1,y_1,x_2,y_2)\longmapsto(u,v),\qquad r_i^2=x_i^2+y_i^2.
\]
Let $J$ be the Jacobian matrix of $\pi$.
Since $\operatorname{rank}J=\operatorname{rank}(JJ^{\mathsf T})$,
direct differentiation gives Table~\ref{tab:almost-toric-local-models}.
All points outside the two listed loci have rank two.

\begin{center}
\begin{minipage}{\linewidth}
\centering
\captionsetup{hypcap=false}
\small
\setlength{\tabcolsep}{5pt}
\renewcommand{\arraystretch}{1.3}
\begin{tabular}{@{}lccc@{}}
\toprule
Local model & $JJ^{\mathsf T}$ & Rank-zero locus & Rank-one locus \\
\midrule
Regular & $I_2$ & $\varnothing$ & $\varnothing$ \\
Rank-one elliptic & $\operatorname{diag}(1,4r_2^2)$
 & $\varnothing$ & $\{r_2=0\}$ \\
Elliptic--elliptic & $\operatorname{diag}(4r_1^2,4r_2^2)$
 & $\{0\}$ & $\{r_1r_2=0\}\setminus\{0\}$ \\
Focus--focus & $(r_1^2+r_2^2)I_2$ & $\{0\}$ & $\varnothing$ \\
\bottomrule
\end{tabular}
\captionof{table}{Ranks of the almost-toric local models.}
\label{tab:almost-toric-local-models}
\end{minipage}
\end{center}

In Table~\ref{tab:almost-toric-local-images}, images and boundaries are local
at $0$: $\R^2$ denotes a full neighborhood, and outer coordinate
boundaries are omitted.

\begin{center}
\begin{minipage}{\linewidth}
\centering
\captionsetup{hypcap=false}
\small
\setlength{\tabcolsep}{5pt}
\renewcommand{\arraystretch}{1.3}
\begin{tabular}{@{}lccc@{}}
\toprule
Local model & Image near $0$ & Boundary near $0$ & Boundary preimage \\
\midrule
Regular & $\R^2$ & $\varnothing$ & $\varnothing$ \\
Rank-one elliptic & $\R\times[0,\infty)$ & $\{v=0\}$ & $\{r_2=0\}$ \\
Elliptic--elliptic & $[0,\infty)^2$ & $\{uv=0\}$ & $\{r_1r_2=0\}$ \\
Focus--focus & $\R^2$ & $\varnothing$ & $\varnothing$ \\
\bottomrule
\end{tabular}
\captionof{table}{Images and boundary preimages of the local models.}
\label{tab:almost-toric-local-images}
\end{minipage}
\end{center}

\begin{example}\label{ex:noncompact-elliptic}
	On standard symplectic $\C^2$, the map $\pi:\C^2\to[0,\infty)^2$,
	$\pi(z_1,z_2)=(|z_1|^2/2,|z_2|^2/2)$, is a proper toric fibration
	\cite[Example~3.11]{Symington2003}. The open edges and corner of the base
	correspond to rank-one elliptic and elliptic--elliptic singularities, respectively.
\end{example}

To our knowledge, the following regularity criterion has not been recorded
in this generality.

\begin{proposition}[Regularity of almost-toric fibrations]\label{prop:almost-toric-regularity}
Let $(M^4,\omega_M=d\lambda_M)$ be a connected exact symplectic manifold
without boundary, homotopy equivalent to $\T^2$.
Every almost-toric fibration $\pi:M\to B$ must be a smooth Lagrangian torus
fibration in the sense of Definition~\ref{def:smooth-fibration}.
\end{proposition}

We prepare the proof with three lemmas. For these lemmas, fix an
almost-toric fibration $\pi:(M^4,\omega_M)\to B$ as in
Definition~\ref{def:almost-toric-fibration}; any additional hypotheses will be stated
explicitly. Homology and cohomology below use $\mathbb F_2$ coefficients.
We set $
E=\pi^{-1}(\partial B)$, $N=M\setminus E$, and $B^\circ=\Int B$.
Write
$P_r=\{p\in M:\operatorname{rank}(d\pi_p)=r\}$ for $r=0,1$.
Thus the critical set is $P_0\cup P_1$, and all remaining points are regular.
Set \(\Delta=\pi(P_0\cap N)\subset B^\circ\), the set of nodal values.
We use this notation in the lemmas and in the final proof of
Proposition~\ref{prop:almost-toric-regularity}.

\begin{lemma}[Structure of the elliptic locus]\label{lem:elliptic-locus}
The set $P_0$ is closed and discrete, and
$P_0\cap\pi^{-1}(C)$ is finite for every compact $C\subset B$.
The elliptic locus $E=\overline{P_1}$ is a closed, locally finite
symplectic normal-crossing divisor. Its crossings are ordinary double
points (nodes of the divisor), occurring precisely at the
elliptic--elliptic points $E\cap P_0$.
\end{lemma}

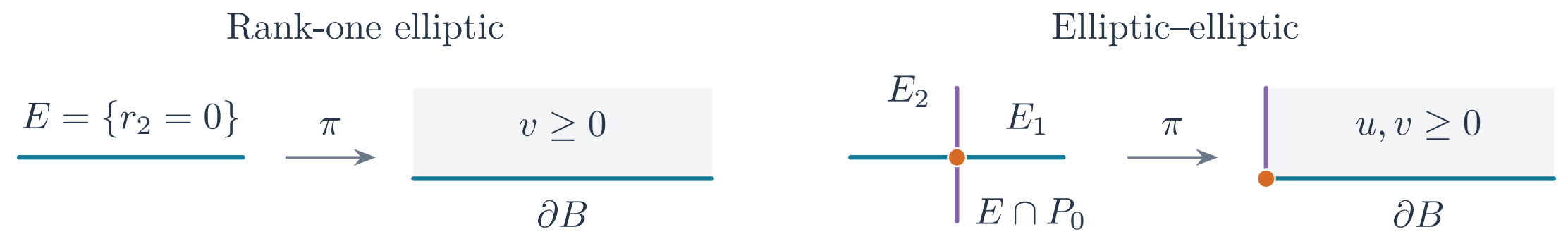
\begin{figure}[htbp]
\centering
\begingroup
\definecolor{elink}{HTML}{26354A}
\definecolor{elmuted}{HTML}{697789}
\definecolor{elblue}{HTML}{167D9A}
\definecolor{elpurple}{HTML}{8765AD}
\definecolor{elpoint}{HTML}{D76B25}
\begin{tikzpicture}[x=1cm,y=1cm,>=Stealth,
  every node/.style={font=\small,text=elink},
  branch/.style={line width=1.2pt,line cap=round},
  projection/.style={->,draw=elmuted,line width=.65pt},
  critical/.style={circle,fill=elpoint,draw=white,line width=.5pt,inner sep=1.7pt}]
\path[use as bounding box] (0,-.15) rectangle (14.8,2.05);

\node at (3.55,1.85) {Rank-one elliptic};
\draw[branch,elblue] (.30,.65) -- (2.40,.65);
\node at (1.35,1.01) {$E=\{r_2=0\}$};
\draw[projection] (2.80,.65) -- (3.65,.65);
\node at (3.22,.94) {$\pi$};
\fill[elmuted!8] (4.00,.45) rectangle (6.80,1.30);
\draw[branch,elblue] (4.00,.45) -- (6.80,.45);
\node at (5.40,.94) {$v\geq0$};
\node at (5.40,.12) {$\partial B$};

\node at (11.15,1.85) {Elliptic--elliptic};
\draw[branch,elblue] (8.10,.65) -- (10.10,.65);
\draw[branch,elpurple] (9.10,.05) -- (9.10,1.30);
\node[anchor=south east] at (10.10,.73) {$E_1$};
\node[anchor=east] at (8.99,1.27) {$E_2$};
\node[critical] at (9.10,.65) {};
\node at (9.79,.12) {$E\cap P_0$};
\draw[projection] (10.70,.65) -- (11.55,.65);
\node at (11.12,.94) {$\pi$};
\fill[elmuted!8] (12.00,.45) rectangle (14.70,1.30);
\draw[branch,elblue] (12.00,.45) -- (14.70,.45);
\draw[branch,elpurple] (12.00,.45) -- (12.00,1.30);
\node[critical] at (12.00,.45) {};
\node at (13.43,.94) {$u,v\geq0$};
\node at (13.43,.12) {$\partial B$};
\end{tikzpicture}
\endgroup
\captionsetup{font=small,skip=3pt}
\caption{Elliptic local models. Colored source lines represent symplectic
surfaces in $M^4$ whose images are boundary edges of $B$.
In the right panel, the branches $E_1=\{r_2=0\}$ and $E_2=\{r_1=0\}$
map to the nonnegative $u$- and $v$-axes, respectively.
Their transverse intersection maps to the corner.}
\label{fig:elliptic-locus}
\end{figure}

\begin{proof}
The set $P_0=\{p:d\pi_p=0\}$ is closed by continuity of $d\pi$.
It is discrete by Table~\ref{tab:almost-toric-local-models}.
Properness makes $\pi^{-1}(C)$ compact for compact $C\subset B$, so
its closed discrete subset $P_0\cap\pi^{-1}(C)$ is finite.

The set $E$ is closed. By Table~\ref{tab:almost-toric-local-images}, regular
and focus--focus points lie over interior points of $B$; conversely,
Definition~\ref{def:almost-toric-fibration} permits only focus--focus critical
points in interior singular fibers. Thus $E$ is exactly the elliptic
locus. The elliptic rows of Table~\ref{tab:almost-toric-local-images} and the
symplectic restrictions above describe it locally as one symplectic
plane or two transverse symplectic planes (see Figure~\ref{fig:elliptic-locus}), proving
the locally finite normal-crossing assertion. At an elliptic--elliptic
point, the coordinates $z_j=x_j+iy_j$ identify $E$ locally with
$\{z_1z_2=0\}$, the standard node. These nodes of $E$ are distinct from
the focus--focus points in $N$. Every other point of $E$ lies in $P_1$.
Each crossing is approached along its rank-one branches, so
$E=\overline{P_1}$.
\end{proof}

\Needspace{6\baselineskip}
\begin{lemma}[Homology of the elliptic complement]\label{lem:elliptic-complement-homology}
If $\omega_M=d\lambda_M$ and $H_3(M)=0$, then every connected component
of $E$ is noncompact and $H_3(N)=0$.
\end{lemma}

\begin{proof}
If $C$ were a compact connected component of $E$, its normalization
$\nu:\widetilde C\to C\subset M$, obtained by separating the branches
at each node, would be a compact, possibly disconnected surface without
boundary with area form $\nu^*\omega_M$. Stokes' theorem would give
$0<\int_{\widetilde C}\nu^*\omega_M
=\int_{\widetilde C}d(\nu^*\lambda_M)=0$, a contradiction.
Thus every component of $E$ is noncompact, and
$H_c^0(E)=0$.

Since $H_c^1(M)\cong H_3(M)=0$ by Poincar\'e duality
\cite[Theorem~3.35]{Hatcher}, the compact-support exact sequence
$H_c^0(E)\to H_c^1(N)\to H_c^1(M)$ and duality on $N$ give
$H_3(N)\cong H_c^1(N)=0$.
\end{proof}

\begin{lemma}[Topology of the base and regular fibers]\label{lem:almost-toric-fiber-topology}
Suppose that $\omega_M$ is exact and $M\simeq\T^2$.
Then $B^\circ\cong\R^2$, and for every regular fiber $F_b=\pi^{-1}(b)$,
the inclusion $F_b\hookrightarrow M$ induces an isomorphism on
fundamental groups.
\end{lemma}

\begin{proof}
Since $M\simeq\T^2$, we have $H_3(M)=0$, and
Lemma~\ref{lem:elliptic-complement-homology} and Poincar\'e duality give
$H_c^1(N)\cong H_3(N)=0$.
Put $p=\pi|_N:N\to B^\circ$. Properness and connected fibers imply
$p_*\underline{\mathbb F}_2=\underline{\mathbb F}_2$.
The compactly supported Leray spectral sequence
\cite[\S\,V.1.6, pp.~49--50]{Borel1984},
\[
E_2^{a,b}=H_c^a(B^\circ;R^bp_*\underline{\mathbb F}_2)
\Longrightarrow H_c^{a+b}(N;\mathbb F_2),
\]
has no differentials entering or leaving $E_2^{1,0}$, so its edge map
gives an injection
$H_c^1(B^\circ)\hookrightarrow H_c^1(N)=0$.
Exactness rules out compact $M$ by Stokes' theorem, and properness
thus rules out compact $B$. Then $B^\circ$ is noncompact, and duality gives $H_1(B^\circ;\mathbb F_2)=0$.
This vanishing rules out handles, crosscaps, and more than one end.
The classification of noncompact surfaces
\cite{Richards1963} then gives $B^\circ\cong\R^2$.

Fix a regular value $b$ and choose a diffeomorphism
$\phi:B^\circ\to\R^2$ with $\phi(b)=0$. The function
\[
h=\|\phi\circ p\|^2:N\longrightarrow[0,\infty)
\]
has compact sublevel sets and minimum set $F_b$.
By Lemma~\ref{lem:elliptic-locus}, every point of $N=M\setminus E$ is either regular or
focus--focus. Since
$d\|\phi\|^2$ is nonzero away from $b$, the critical points of $h$
outside $F_b$ are precisely the focus--focus points. At each such
point, the local model shows that the quadratic part of $h$ is
$a(x_1y_1+x_2y_2)+c(x_1y_2-x_2y_1)$ with $(a,c)\ne(0,0)$.
Its Hessian has signature $(2,2)$, so every critical point outside
$F_b$ has index two. The relative Morse-theoretic lemma
\cite[Lemma~22.5, p.~121]{Milnor1963} gives
$\pi_1(N,F_b)=0$. Thus, the homotopy exact sequence of the pair \((N,F_b)\) shows that the inclusion induces a surjection \(\pi_1(F_b)\twoheadrightarrow\pi_1(N)\).

Fix $x\in F_b\subset N$. Since $E$ is locally finite with strata of
codimension at least two, general position lets us perturb every loop
based at $x$ off $E$, fixing $x$. The perturbed
loop therefore lies in $N$, proving that $\pi_1(N,x)\to\pi_1(M,x)$
is surjective. Hence the composite $\pi_1(F_b,x)\to\pi_1(M,x)$ is
surjective as well. Since $F_b\cong\T^2$ and $M\simeq\T^2$, both
groups are isomorphic to $\Z^2$. Since every surjective homomorphism
$\Z^2\to\Z^2$ is injective, the inclusion induces an isomorphism as desired.
\end{proof}

\begin{proof}[Proof of Proposition~\ref{prop:almost-toric-regularity}]
By Lemma~\ref{lem:almost-toric-fiber-topology}, the inclusion of each regular
fiber into $M$ induces an injection on fundamental groups. We show
that a nodal or elliptic fiber would contradict this injectivity.

If a nodal fiber exists, choose a focus--focus point $p$ on it.
In Darboux coordinates $(x_1,y_1,x_2,y_2)$ centered at $p$ and base
coordinates $(u,v)$ centered at $\pi(p)$, the model in
Definition~\ref{def:almost-toric-fibration} gives
$(u,v)\circ\pi=(x_1y_1+x_2y_2,\ x_1y_2-x_2y_1)$.
On the disk
$
D_\varepsilon=\{y_1=x_1,\ y_2=x_2,\ x_1^2+x_2^2\le\varepsilon\}$
this reduces to $(x_1^2+x_2^2,0)$, so $\partial D_\varepsilon$ lies in
the fiber $F$ over $(u,v)=(\varepsilon,0)$. For sufficiently small
$\varepsilon>0$ this fiber is regular by Lemma~\ref{lem:elliptic-locus}.
The Hamiltonian flow of $v\circ\pi$ rotates $(x_1,x_2)$ and $(y_1,y_2)$
together, so its orbit through $(x_1,x_2)=(y_1,y_2)=(\sqrt\varepsilon,0)$
traces the vanishing cycle $\partial D_\varepsilon$ once per period. This cycle represents
a nontrivial element of $\pi_1(F)$ as in
\cite[\S\,3, paragraph immediately before Lemma~1]{Zung1997}.
Since it bounds $D_\varepsilon$ in $M$, this contradicts the injectivity
of $\pi_1(F)\to\pi_1(M)$ in Lemma~\ref{lem:almost-toric-fiber-topology}.
Thus $\Delta=\varnothing$.

To exclude elliptic fibers, let $C=\pi^{-1}(b_\partial)$ be a fiber over
a boundary value. By Definition~\ref{def:almost-toric-fibration}, $C$ is a circle
over an edge or a point over a corner. Choose a tubular neighborhood $U$
of $C$ in $M$. It deformation retracts onto $C$, so $\pi_1(U)\cong\pi_1(C)$
is either $\Z$ or trivial.
Since a proper continuous map between manifolds is a closed map, the set \(\pi(M\setminus U)\) is closed in \(B\).
This set omits
$b_\partial$ as its entire fiber $C$ lies in $U$. Thus $F_b\subset U$
for every $b$ sufficiently close to $b_\partial$. Choose such an interior value $b$.
Since $\Delta=\varnothing$, its fiber $F_b$ is a regular torus.
The map $\pi_1(F_b)\cong\Z^2\to\pi_1(M)$ factors through $\pi_1(U)$
and hence cannot be injective, contradicting Lemma~\ref{lem:almost-toric-fiber-topology}.
Thus $\partial B=\varnothing$ and
$\pi$ is a smooth Lagrangian torus fibration.
\end{proof}

\section{Lagrangian shape and smooth torus fibrations}\label{sec:shape}

Having reduced the almost-toric case to the smooth case in
Section~\ref{sec:dimension-four}, we now establish the shape--volume identity
for smooth Lagrangian torus fibrations.
The results of this section hold in arbitrary dimension.

Throughout this section, let $n\geq1$ and let
$(M^{2n},\omega_M=d\lambda_M)$ be a connected exact symplectic manifold
with $H_1(M;\Z)\cong\Z^n$.
The one-form $\lambda_M$ is a chosen primitive of $\omega_M$.
Examples include \(T^*\mathbb T^n\) and its product domains \(D\times\mathbb T^n\), with the canonical Liouville form, where \(D\subset\mathbb R^n\) is nonempty, open, and contractible.

\begin{definition}\label{def:shape}
A \emph{marking} of $M$ is an ordered integral basis of $H^1(M;\Z)$.
Fixing a marking identifies $H^1(M;\R)$ with $\R^n$, with coordinates
$(p_1,\ldots,p_n)$ and the associated Lebesgue measure $\Vol_n$.
An embedded Lagrangian $n$-torus $\iota:L\hookrightarrow(M,\omega_M)$ is \emph{essential} if
$\iota_*:H_1(L;\Z)\to H_1(M;\Z)$ is an isomorphism.
Since $d(\iota^*\lambda_M)=\iota^*\omega_M=0$, its action class with respect to $\lambda_M$ is
\[
I_{\lambda_M}(L)=(\iota^*)^{-1}[\iota^*\lambda_M]\in H^1(M;\R).
\]
The components $p_1,\ldots,p_n$ of $I_{\lambda_M}(L)$ in this marking are called
the \emph{marked action coordinates} of $L$.
The \emph{marked essential Lagrangian shape} is
\[
\Shape_{\lambda_M}(M)
=\{I_{\lambda_M}(L):L\subset M\text{ is an essential Lagrangian }n\text{-torus}\}.
\]
\end{definition}

This is the Sikorav--Eliashberg Lagrangian shape
\cite{Sikorav1989,Sikorav1991,Eliashberg}, associated with a fixed isomorphism on
integral first cohomology and transported to $H^1(M;\R)$; see also the work of M\"uller
\cite[Definitions~2.1 and~2.8]{Muller}. Here ``essential'' emphasizes the
homological condition in Definition~\ref{def:shape}.

The shape is open by \cite[Lemma~2.4]{Muller}, and hence measurable.
If $\beta$ is a closed one-form on $M$, changing the primitive translates
the action classes and the shape \cite[Lemma~2.2]{Muller}:
\begin{equation}\label{eq:primitive-shift}
I_{\lambda_M+\beta}(L)=I_{\lambda_M}(L)+[\beta],
\qquad
\Shape_{\lambda_M+\beta}(M)=\Shape_{\lambda_M}(M)+[\beta].
\end{equation}
Changing the integral marking transforms the coordinates of the shape by a matrix
$A\in\mathrm{GL}(n,\Z)$; see \cite[Proposition~2.14]{Muller}.
Since $\det A=\pm1$, this preserves its volume as $\Vol_n(A\Shape)
=|\det A|\,\Vol_n(\Shape)
=\Vol_n(\Shape)$.
Changing the primitive only translates the shape by \textup{(\ref{eq:primitive-shift})},
so its $n$-dimensional volume is independent of both the chosen primitive
and the integral marking.

Consider the cotangent bundle $T^*\T^n=\R^n\times\T^n$, with coordinates
$(a,\theta)$, where $\theta_j\in\R/\Z$ for $j=1,\ldots,n$. Its canonical
Liouville form and symplectic form are
\begin{equation}\label{eq:can-form}
\lambda_{\mathrm{can}}=\sum_{j=1}^n a_j\,d\theta_j,
\qquad
\omega_{\mathrm{can}}=d\lambda_{\mathrm{can}}
=\sum_{j=1}^n da_j\wedge d\theta_j.
\end{equation}
We use its standard marking $[d\theta_1],\ldots,[d\theta_n]$.
For $c=(c_1,\ldots,c_n)\in\R^n$, we define the translation
$\tau_c(a,\theta):=(a+c,\theta)$.
Then
\begin{equation}\label{eq:log-transformations}
\tau_c^*\lambda_{\mathrm{can}}
=\lambda_{\mathrm{can}}+\sum_{j=1}^n c_j\,d\theta_j.
\end{equation}

To prove injectivity of the action map in Proposition~\ref{prop:shape-volume}, we will use
Gromov's exact-isotopy intersection theorem under a bounded geometry hypothesis
\cite[\S\,2.3.B$'_3$, pp.~330--331]{Gromov}.

\begin{theorem}[Gromov]\label{thm:exact-intersection}
Let $(N,\omega_N)$ be a symplectic manifold admitting an $\omega_N$-compatible
almost complex structure $J$ such that the metric
$g=\omega_N(\cdot,J\cdot)$ is complete, has bounded sectional curvature,
and has injectivity radius bounded away from zero. Let
$\iota_t:L\hookrightarrow(N,\omega_N)$, $0\leq t\leq1$, be an exact Lagrangian
isotopy of a nonempty closed manifold. If every $\iota_t(L)$ is weakly exact, then
\[
\iota_0(L)\cap\iota_1(L)\ne\varnothing.
\]
\end{theorem}

Here the isotopy is \emph{exact} if the one-forms
$\alpha_t(v)=\omega_N\bigl(\partial_t\iota_t,d\iota_t(v)\bigr)$, $v\in TL$,
are exact on $L$. The image $\iota_t(L)$ is \emph{weakly exact} if
$\int_{D^2}u^*\omega_N=0$ for every smooth map
$u:(D^2,\partial D^2)\to(N,\iota_t(L))$.

We apply this theorem in $(T^*\T^n,\omega_{\mathrm{can}})$, whose flat compatible
metric satisfies the geometric hypotheses. If each
$\iota_t^*\lambda_{\mathrm{can}}$ is exact, Cartan's formula gives
$[\alpha_t]=\frac{d}{dt}[\iota_t^*\lambda_{\mathrm{can}}]=0$, and Stokes' theorem
gives weak exactness.

\begin{lemma}[Shape of product domains]\label{lem:log-domain-shape}
Let $D\subset\R^n$ be a nonempty contractible open set and put
$U_D=D\times\T^n\subset T^*\T^n$. In the standard marking,
we have $\Shape_{\lambda_{\mathrm{can}}}(U_D)=D$.
\end{lemma}

\begin{proof}
This is essentially Sikorav’s shape theorem; see \cite[Theorem~2.17]{Muller}. We give a brief proof for completeness.
	Every coordinate torus $\{a=p\}$, $p\in D$, is essential and has action $p$, so
	$D\subset\Shape_{\lambda_{\mathrm{can}}}(U_D)$.
	Conversely, suppose for contradiction that there exists an essential Lagrangian torus \(L\subset U_D\) whose action vector \(p=I_{\lambda_{\mathrm{can}}}(L)\) lies outside \(D\).
	By \textup{(\ref{eq:log-transformations})}, $K=\tau_{-p}(L)$ is a closed exact
	Lagrangian in $T^*\T^n$. It is disjoint from the zero section $\{a=0\}$, because
	$K\subset U_{D-p}$ and $0\notin D-p$. Compactness gives constants
	$0<m\leq R$ such that
	$m\leq|a|\leq R$ on $K$.
	For $s>0$, let $\delta_s:T^*\T^n\to T^*\T^n$ be the fiberwise dilation
	$\delta_s(a,\theta)=(sa,\theta)$. If $r>R/m$, then $|a|\geq rm>R$
	on $\delta_r(K)$, so $\delta_r(K)\cap K=\varnothing$.
	The family $\delta_\sigma(K)$, $1\leq\sigma\leq r$, is a Lagrangian isotopy,
	and each member is exact since
	$\delta_\sigma^*\lambda_{\mathrm{can}}=\sigma\lambda_{\mathrm{can}}$.
	The endpoints are disjoint, contradicting Theorem~\ref{thm:exact-intersection}.
\end{proof}

\begin{proposition}[Shape-volume identity]\label{prop:shape-volume}
Let \((M^{2n},\omega_M=d\lambda_M)\) be as above.
Assume that $M$ is homotopy equivalent to $\T^n$ and that there is an exact
symplectic embedding
$j:(M,\omega_M)\hookrightarrow(T^*\T^n,\omega_{\mathrm{can}})$
inducing an isomorphism on integral first homology. Here exactness means that
$j^*\lambda_{\mathrm{can}}=\lambda_M+dh$ for a smooth function $h$ on $M$.
If $(M,\omega_M)$ admits a smooth Lagrangian torus fibration in the sense of Definition~\ref{def:smooth-fibration}, then
\begin{equation}\label{eq:shape-volume}
\Vol_n\Shape_{\lambda_M}(M)=\int_M\frac{\omega_M^n}{n!}.
\end{equation}
\end{proposition}

\begin{proof}
Let $\pi:M\to B$ be such a fibration as in Definition~\ref{def:smooth-fibration}.
Its fibers are diffeomorphic to $\T^n$ by the Arnold--Liouville theorem, and
Ehresmann's theorem makes $\pi$ a smooth torus bundle. We first prove that $B$ is
contractible and that every fiber is essential.
Fix a fiber $F\simeq\T^n$. Since $M\simeq\T^n$, both $F$ and $M$ are aspherical.
The homotopy sequence of the bundle therefore gives $\pi_k(B)=0$ for $k\geq3$
and an exact sequence
\[
0\longrightarrow\pi_2(B)\longrightarrow\pi_1(F)
\longrightarrow\pi_1(M)
\longrightarrow\pi_1(B)\longrightarrow0.
\]
In particular, $\pi_2(B)$ is a subgroup of $\Z^n$, hence isomorphic to $\Z^r$
for some $0\leq r\leq n$. The universal cover $\widetilde B$ is a $K(\Z^r,2)$, hence
$\widetilde B\simeq(\mathbb{CP}^{\infty})^r$
by \cite[Proposition~4.30 and Example~4.50]{Hatcher}.
If $r>0$, this space has nonzero cohomology in arbitrarily high degrees
\cite[Theorem~3.19]{Hatcher}, contradicting the vanishing of cohomology
above degree $n$ for the smooth $n$-manifold $\widetilde B$.
Thus $r=0$, and $\widetilde B$ is contractible.
The displayed sequence shows that $\pi_1(B)$ is finite. As $B$ is an
aspherical manifold, $\pi_1(B)$ is also torsion-free by
\cite[Proposition~2.45]{Hatcher}, hence trivial. Thus $B$ is contractible, and
the same sequence shows that every fiber inclusion induces isomorphisms on
fundamental groups and integral first homology.

The volume of shape is independent of the integral marking. We use the basis
$j^*[d\theta_1],\ldots,j^*[d\theta_n]$ on $M$ and transport it to every fiber by
the fiber inclusions. Define the action map
\[
A:B\longrightarrow\R^n,
\qquad A(\xi)=I_{\lambda_M}(\pi^{-1}(\xi)).
\]
The map is smooth. Local action--angle coordinates show that $dA$ is invertible
everywhere. We now prove that $A$ is injective:
Suppose $\xi_0\ne\xi_1$ and $A(\xi_0)=A(\xi_1)$. Choose a smooth path
$\xi(t)$, $0\leq t\leq1$, from $\xi_0$ to $\xi_1$ and a smooth family of
parametrizations $f_t:\T^n\to\pi^{-1}(\xi(t))$ respecting the fixed integral
marking. Write $p(t)=A(\xi(t))$ and define
\[
\iota_t=\tau_{-p(t)}\circ j\circ f_t,
\qquad K_t=\iota_t(\T^n).
\]
Each $\iota_t^*\lambda_{\mathrm{can}}$ is exact: the identity
$j^*\lambda_{\mathrm{can}}=\lambda_M+dh$ identifies the action classes, and
\textup{(\ref{eq:log-transformations})} shows that $\tau_{-p(t)}$ subtracts precisely
$p(t)$ in the chosen marking. Thus Theorem~\ref{thm:exact-intersection} implies
$K_0\cap K_1\ne\varnothing$. But the original fibers are disjoint, $j$ is an
embedding, and the endpoint translations are identical, so $K_0$ and $K_1$ are
disjoint. This contradiction proves injectivity.

The injective local diffeomorphism $A$ identifies $B$ with a contractible
open set $D\subset\R^n$. By the global action--angle theorem
\cite{Duistermaat}, there is a symplectomorphism
$\Psi:(U_D,\omega_{\mathrm{can}}) \to (M,\omega_M)$ where $U_D=D\times\T^n$. Thus Lemma~\ref{lem:log-domain-shape} gives
\[
\Vol_n\Shape_{\lambda_M}(M)
=\Vol_n\Shape_{\lambda_{\mathrm{can}}}(U_D)
=\Vol_n(D)=\int_{U_D}\frac{\omega_{\mathrm{can}}^n}{n!}
=\int_M\frac{\omega_M^n}{n!}.\qedhere
\]
\end{proof}

\section{The symplectic obstruction}\label{sec:tori}

Using the logarithmic coordinates of Section~\ref{sec:geometry}, we regard
$\lambda_{\mathrm{can}}$ and $\omega_{\mathrm{can}}$ from
\textup{(\ref{eq:can-form})} as forms on $X$. All action classes on $X$ are
expressed in the angular marking $[d\theta_1],[d\theta_2]$ of $H^1(X;\Z)$.
The map
\[
\rho:\C\longrightarrow\C,
\qquad \rho(z)=\frac{z}{\sqrt{1+|z|^2}}
\]
is odd and is a smooth
embedding.

The standard Euclidean primitive on $\C^2$ is
\begin{equation}\label{eq:flat-primitive}
\lambda_0=d^c(|z_1|^2+|z_2|^2)
=\frac{1}{2\pi}\operatorname{Im}(\bar z_1\,dz_1+\bar z_2\,dz_2).
\end{equation}
A direct calculation gives
\begin{equation}\label{eq:radial-exact}
(\rho\times\rho)^*\lambda_0
=d^c\bigl(\log(1+|z_1|^2)+\log(1+|z_2|^2)\bigr).
\end{equation}
In logarithmic coordinates on $X$, we have
$\lambda_0=e^{2a_1}\,d\theta_1+e^{2a_2}\,d\theta_2$.
Define the radial map in these coordinates by
\begin{equation}\label{eq:flat-to-log}
H:X\longrightarrow X,
\qquad
H(a_1,a_2,\theta_1,\theta_2)
= (e^{2a_1},e^{2a_2},\theta_1,\theta_2 ).
\end{equation}
Then $H^*\lambda_{\mathrm{can}}=\lambda_0$, and $H$ is a diffeomorphism onto
$\{(a,\theta)\in X:a_1>0,\ a_2>0\}$. Hence
$H:(X,\lambda_0)\hookrightarrow(X,\lambda_{\mathrm{can}})$ is an exact symplectic
embedding preserving the angular coordinates.

\Needspace{7\baselineskip}
\begin{proposition}
	\label{prop:log-embedding}
There is an exact symplectic embedding
$j:(X,\lambda)\hookrightarrow(X,\lambda_{\mathrm{can}})$
whose induced map $j_*:H_1(X;\Z)\to H_1(X;\Z)$ is the identity.
In fact, $j^*\lambda_{\mathrm{can}}=\lambda$.
\end{proposition}

\begin{proof}
We first construct an exact symplectic embedding from $(X,\lambda)$ into
$(X,\lambda_0)$; composing it with $H$ will then give the desired map.
Set $F=T^{-1}\circ(\rho\times\rho)\circ T$, where $T=(u,v)$ is given by \textup{(\ref{eq:unitary})}. Since $\rho\times\rho$ is a smooth
embedding and $T,T^{-1}$ are diffeomorphisms, $F$ is a smooth embedding of $\C^2$
into itself. 
Since $\lambda_0=d^c(|z_1|^2+|z_2|^2)$ and holomorphic pullback commutes
with $d^c$, unitarity of $T$ and $T^{-1}$ gives
$T^*\lambda_0=(T^{-1})^*\lambda_0=\lambda_0$ as identities of forms on \(\mathbb C^2\).
Put $G(w)=\sum_{j=1}^2\log(1+|w_j|^2)$. By \textup{(\ref{eq:radial-exact})},
$(\rho\times\rho)^*\lambda_0=d^cG$, and \textup{(\ref{eq:factorization})} gives
$G\circ T=\Phi$. Therefore,
\[
F^*\lambda_0
=T^*\bigl((\rho\times\rho)^*\lambda_0\bigr)=T^*(d^cG) =d^c(G\circ T)=d^c\Phi=\lambda.
\]

To compose $F:\mathbb C^2\to\mathbb C^2$ with $H$, we must also check that $F(X)\subset X$. Write $F=(F_1,F_2)$. The formula
for $T^{-1}$ gives
\[
F_1=\frac{\rho(u)+\rho(v)}{\sqrt2},
\qquad
F_2=\frac{\rho(u)-\rho(v)}{i\sqrt2}.
\]
Oddness and injectivity of $\rho$ therefore imply
\[
\begin{aligned}
F_1=0
&\iff \rho(u)=\rho(-v)\iff u=-v\iff z_1=0,\\
F_2=0
&\iff \rho(u)=\rho(v)\iff u=v\iff z_2=0.
\end{aligned}
\]
Thus neither coordinate of $F(z_1,z_2)$ vanishes when $(z_1,z_2)\in X$.
Consequently, $j=H\circ F|_X:X\to X$ is a smooth embedding, and
\[
j^*\lambda_{\mathrm{can}}=(F|_X)^*(H^*\lambda_{\mathrm{can}})=(F|_X)^*\lambda_0=\lambda.
\]

It remains to check the induced map on first homology. For $\varepsilon\in[0,1]$, define
\[
\rho_\varepsilon(z)=\frac{z}{\sqrt{1+\varepsilon|z|^2}},
\qquad F_\varepsilon=T^{-1}\circ(\rho_\varepsilon\times\rho_\varepsilon)\circ T.
\]
Every $\rho_\varepsilon$ is odd and injective. Hence every $F_\varepsilon$ restricts to a map $X\to X$, with
$F_0=\id$ and $F_1=F$. This homotopy proves that $F$ induces the identity on
$H_1(X;\Z)$. Since $H$ preserves the angles, $j$ has the same property.
\end{proof}

Note that Proposition~\ref{prop:log-embedding} supplies the embedding hypothesis of Proposition~\ref{prop:shape-volume}.
We construct contractions of the shape toward points of the open set
$\Delta^\circ=\{(\alpha,\beta):\alpha>0,\ \beta>0,\ \alpha+\beta<1\}$.

\begin{proposition}[Contraction centers for the example]\label{prop:contraction-centers}
For $\kappa\in\Delta^\circ$, $v\in\overline{\Shape_\lambda(X)}$, and
$0<\sigma<1$, one has
$(1-\sigma)\kappa+\sigma v\in\Shape_\lambda(X)$.
\end{proposition}

\begin{proof}
Fix $\kappa=(\alpha,\beta)\in\Delta^\circ$ and set
$\psi_\kappa=\Phi-\alpha\log|z_1|^2-\beta\log|z_2|^2$.
Since $d^c\log|z_j|^2=d\theta_j$ and $dd^c\log|z_j|^2=0$ on $X$,
we have $dd^c\psi_\kappa=\omega$ and $
\lambda_\kappa:=d^c\psi_\kappa
=\lambda-\alpha\,d\theta_1-\beta\,d\theta_2$.
Put $s=|z_1|^2$, $t=|z_2|^2$, $\gamma=1-\alpha-\beta>0$, and
$S=e^\Phi\geq1+s+t$. Then
\[
\psi_\kappa
=\gamma\log S+\alpha\log(S/s)+\beta\log(S/t).
\]
All three terms are nonnegative. On each sublevel set of $\psi_\kappa$,
$S,S/s,S/t$ are bounded, so $s,t$ are bounded above and away from zero.
By continuity, the sublevel set is therefore compact in $X$.
Let $Z_\kappa$ be the Liouville vector field with
$\iota_{Z_\kappa}\omega=\lambda_\kappa$.
For the standard complex structure $J$, our convention
$d^cf=-(4\pi)^{-1}df\circ J$ gives
$d\psi_\kappa(-Z_\kappa)=-4\pi\,\omega(Z_\kappa,JZ_\kappa)\leq0$.
Thus negative trajectories stay in compact sublevel sets, so the flow of
$-Z_\kappa$ exists for every $\tau\geq0$.
Its time-$\tau$ map $F^\kappa_\tau:X\hookrightarrow X$ satisfies the
Liouville flow identities \cite[\S\,2.8, equation~(2.16)]{CieliebakEliashberg2012}:
\begin{equation}\label{eq:flow-scaling}
(F^\kappa_\tau)^*\lambda_\kappa=e^{-\tau}\lambda_\kappa,
\qquad (F^\kappa_\tau)^*\omega=e^{-\tau}\omega.
\end{equation}
Each $F^\kappa_\tau$ is homotopic to the identity and takes essential
Lagrangian tori to essential Lagrangian tori, preserving the angular marking.
For such a torus $L$ with $p=I_\lambda(L)$, these identities and
\textup{(\ref{eq:primitive-shift})} thus give
$
I_\lambda(F^\kappa_\tau(L))
=\kappa+e^{-\tau}I_{\lambda_\kappa}(L)
=\kappa+e^{-\tau}(p-\kappa)$.
By Definition~\ref{def:shape}, taking $\sigma=e^{-\tau}$ proves
\begin{equation}\label{eq:shape-contraction}
(1-\sigma)\kappa+\sigma p\in\Shape_\lambda(X)
\qquad
\bigl(\kappa\in\Delta^\circ,\ p\in\Shape_\lambda(X),\ 0<\sigma\leq1\bigr).
\end{equation}
Finally, fix $v\in\overline{\Shape_\lambda(X)}$ and $0<\sigma<1$.
Since $\Delta^\circ$ is open, choose $p\in\Shape_\lambda(X)$ sufficiently
close to $v$ that $\kappa'=\kappa+\frac{\sigma}{1-\sigma}(v-p)\in\Delta^\circ$.
Applying \textup{(\ref{eq:shape-contraction})} with center $\kappa'$ gives
\[
(1-\sigma)\kappa+\sigma v
=(1-\sigma)\kappa'+\sigma p\in\Shape_\lambda(X).\qedhere
\]
\end{proof}

By Proposition~\ref{prop:explicit-tori}, every torus $L_{c,R}$ is essential in the
sense of Definition~\ref{def:shape}. For $0<c<2$ and $0<R<1$, we compute the action
classes $I_\lambda(L_{c,R})$ and their limits as $R\to0$.
Using the parametrization \textup{(\ref{eq:torus-parameterization})}, define
\[
\gamma_{c,R}(\varphi)=f_{c,R}(0,\varphi),
\qquad J(c,R)=\int_{\gamma_{c,R}}\lambda.
\]
The loop $\gamma_{c,R}$ has angular winding vector $(0,1)$: it winds once in $z_2$
and zero times in $z_1$. Thus $J(c,R)$ is the second marked action coordinate.

\begin{proposition}[Action classes and their limits]\label{prop:torus-actions}
For $0<c<2$ and $0<R<1$, the marked action of $L_{c,R}$ is
\begin{equation}\label{eq:action-J}
I_\lambda(L_{c,R})=(c-J(c,R),J(c,R)).
\end{equation}
For every fixed $c\in(0,2)$, we have $J(c,R)\to0$ as $R\to0$, and hence
$I_\lambda(L_{c,R})\to(c,0)$. In particular,
$(2,0),(0,2)\in\overline{\Shape_\lambda(X)}$.
\end{proposition}

\begin{proof}
The integral of $\lambda$ along the $\theta_1$-coordinate loop equals $c$, since its
tangent vector is $Y$ and $\lambda(Y)=\mu=c$. The integrals of $\lambda$ along the
winding cycles $(1,1)$ and $(0,1)$ are therefore
$c$ and $J(c,R)$, respectively. Hence, if the marked action is $(p_1,p_2)$, then
$p_1+p_2=c$ and $p_2=J(c,R)$, proving \textup{(\ref{eq:action-J})}.

It remains to determine the limit as $R\to0$. Fix $c\in(0,2)$ and recall
$k=1+R^2$ and $h(\varphi)=\tfrac14|1+R^2e^{4\pi i\varphi}|^2$
from the proof of Proposition~\ref{prop:explicit-tori}. Rewriting
\textup{(\ref{eq:radial-moment})} with $\mu=c$ as
$(2-c)hr^2+(1-c)kr-c=0$ and solving for the positive root gives
\[
r_{c,R}(\varphi)
=\frac{2c}{(1-c)k+\sqrt{(1-c)^2k^2+4c(2-c)h(\varphi)}}.
\]
At $R=0$, we have $k=1$ and $h=1/4$, so the expression under the square root
equals $1$ and the denominator equals $2-c>0$. This formula therefore extends
$r_{c,R}(\varphi)$ smoothly to $R=0$, with constant value $r_0=2c/(2-c)$.
Since $\T=\R/\Z$ is compact, $r_{c,R}\to r_0$ and
$\partial_\varphi r_{c,R}\to0$ uniformly in $\varphi$. Consequently,
$\gamma_{c,R} \to  (\sqrt{r_0},0)$ in $C^1(\T,\C^2)$ as $R\to0$.
By \textup{(\ref{eq:factorization})}, $\lambda=d^c\Phi$ extends smoothly to $\C^2$
and is therefore bounded near $(\sqrt{r_0},0)$. The $C^1$ convergence to this
constant loop implies that the Euclidean lengths of $\gamma_{c,R}$ tend to zero.
Hence $J(c,R)=\int_{\gamma_{c,R}}\lambda\to0$.

Thus $(c,0)\in\overline{\Shape_\lambda(X)}$ for every $0<c<2$.
Letting $c\to2$ gives $(2,0)$ in this closure. Interchanging $z_1,z_2$
preserves $\Phi$ and exchanges the marked action coordinates, giving $(0,2)$ as well.
\end{proof}

Set
$T_1=\conv\{(0,0),(2,0),(0,1)\}$ and $
T_2=\conv\{(0,0),(1,0),(0,2)\}$.
The two triangles and their intersection are shown in Figure~\ref{fig:shape} below.

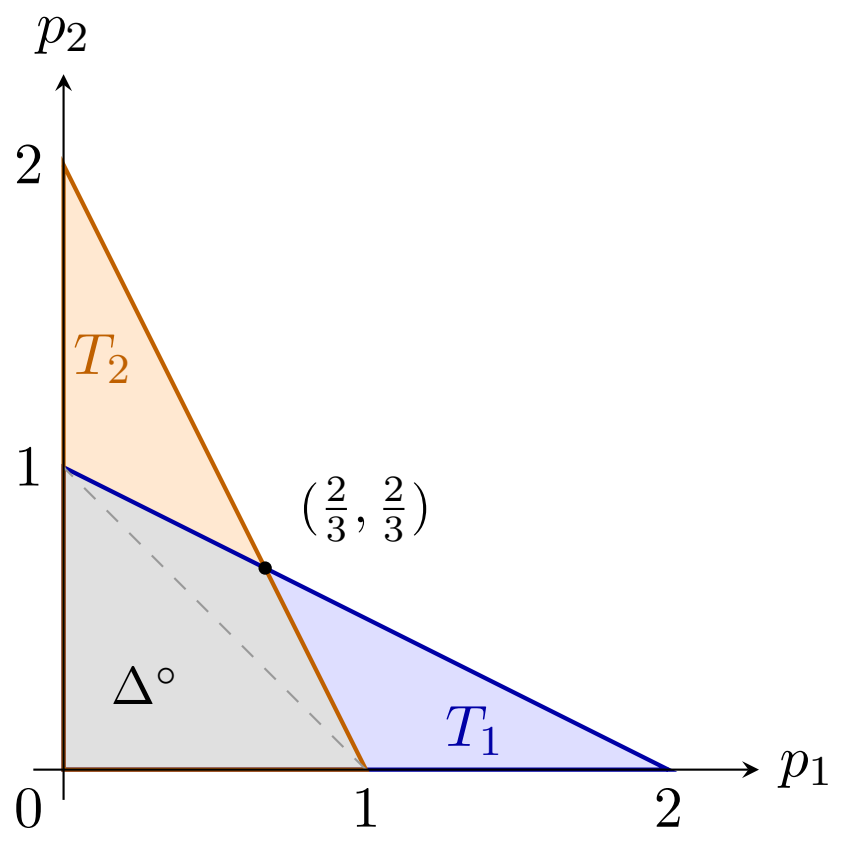
\begin{figure}[htbp]
\centering
\begin{tikzpicture}[scale=2.0,>=stealth]
  \fill[blue!13] (0,0)--(2,0)--(0,1)--cycle;
  \fill[orange!18] (0,0)--(1,0)--(0,2)--cycle;
  \fill[gray!24] (0,0)--(1,0)--(0.666667,0.666667)--(0,1)--cycle;
  \draw[blue!65!black,thick] (0,0)--(2,0)--(0,1)--cycle;
  \draw[orange!75!black,thick] (0,0)--(1,0)--(0,2)--cycle;
  \draw[dashed,gray!80] (1,0)--(0,1);
  \draw[->] (-0.1,0)--(2.3,0) node[right] {$p_1$};
  \draw[->] (0,-0.1)--(0,2.3) node[above] {$p_2$};
  \node[below left] at (0,0) {$0$};
  \node[below] at (1,0) {$1$};
  \node[below] at (2,0) {$2$};
  \node[left] at (0,1) {$1$};
  \node[left] at (0,2) {$2$};
  \fill (0.666667,0.666667) circle (0.022);
  \node[above right,font=\small] at (0.71,0.69) {$(\frac23,\frac23)$};
  \node[blue!65!black] at (1.36,0.13) {$T_1$};
  \node[orange!75!black] at (0.13,1.36) {$T_2$};
  \node[font=\small] at (0.27,0.28) {$\Delta^\circ$};
\end{tikzpicture}
\captionsetup{font=scriptsize}
\caption{The triangles $T_1,T_2$ have union area $4/3$, and their interiors lie in the shape.}
\label{fig:shape}
\end{figure}

\begin{proposition}[The shape estimate]\label{prop:area-bound}
	The marked essential Lagrangian shape satisfies
	\[
	\Int(T_1)\cup\Int(T_2)\subset\Shape_\lambda(X),
	\qquad
	\Area\Shape_\lambda(X)\geq\frac43.
	\]
\end{proposition}

\begin{proof}
The closed convex hull of $\overline{\Delta^\circ}$ and $(2,0)$ is $T_1$. Every point in
$\Int(T_1)$ can be written as $(1-\sigma)\kappa+\sigma(2,0)$ with $\kappa\in\Delta^\circ$ and
$0<\sigma<1$. For instance, extend the segment from $(2,0)$ through the prescribed interior
point until it enters $\Delta^\circ$. By Proposition~\ref{prop:torus-actions},
$(2,0)\in\overline{\Shape_\lambda(X)}$.
Proposition~\ref{prop:contraction-centers}, applied with $v=(2,0)$ and the chosen
$\kappa\in\Delta^\circ$ and $0<\sigma<1$, gives
$(1-\sigma)\kappa+\sigma(2,0)\in\Shape_\lambda(X)$.
Thus $\Int(T_1)\subset\Shape_\lambda(X)$. The same argument with $v=(0,2)$ yields
$\Int(T_2)\subset\Shape_\lambda(X)$.
Now, each triangle has area $1$. Their intersection is the quadrilateral with vertices
$(0,0)$, $(1,0)$, $(2/3,2/3)$, and $(0,1)$; its area is $2/3$. Thus
$
\Area\bigl(\Int(T_1)\cup\Int(T_2)\bigr)
=1+1-\frac23=\frac43$.
\end{proof}

\begin{proof}[Proof of Theorem~\ref{thm:main}]
Suppose that $(X,\omega)$ admits an almost-toric fibration.
By Proposition~\ref{prop:kahler-volume}, $\omega$ is exact and has total symplectic
volume $1$, and the logarithmic coordinates give $X\simeq\T^2$.
Thus Proposition~\ref{prop:almost-toric-regularity} implies that the fibration is
smooth. By Proposition~\ref{prop:log-embedding}, there is an exact symplectic
embedding $j:(X,\omega)\hookrightarrow(T^*\T^2,\omega_{\mathrm{can}})$
inducing an isomorphism on integral first homology. Applying
Proposition~\ref{prop:shape-volume} with $M=X$, $n=2$, and this embedding $j$ gives
$\Area\Shape_\lambda(X)=\int_X\frac{\omega^2}{2}=1$,
contradicting Proposition~\ref{prop:area-bound}.
\end{proof}

\begin{proof}[Proof of Corollary~\ref{cor:higher-dimensional}]
For every $n\ge2$, set $X_n=(\C^*)^n$ and
\[
\Phi_n(z)=\Phi(z_1,z_2)+\sum_{j=3}^n\log(1+|z_j|^2),
\qquad \lambda_n=d^c\Phi_n,\qquad \omega_n=d\lambda_n.
\]
The product potential makes $\omega_n$ real-analytic, exact, and K\"ahler.
Each added factor is exactly symplectomorphic to
$((0,1)\times\T,p\,d\theta)\subset T^*\T$ via
$p=|z|^2/(1+|z|^2)$ and $\theta=\arg z/(2\pi)$, and has area $1$.
Thus Propositions~\ref{prop:kahler-volume} and~\ref{prop:log-embedding} give total volume $1$
and a product exact embedding $X_n\hookrightarrow T^*\T^n$ inducing an
isomorphism on integral first homology. Also $X_n\simeq\T^n$.
Products of essential tori in $X$ with coordinate circles in the added
factors give, in the angular marking,
\[
\Shape_{\lambda_n}(X_n)
\supset\Shape_\lambda(X)\times(0,1)^{n-2},
\qquad
\Vol_n\Shape_{\lambda_n}(X_n)\ge\frac43
\]
by Proposition~\ref{prop:area-bound}. A smooth Lagrangian torus fibration would
instead force this shape volume to equal the symplectic volume $1$ by
Proposition~\ref{prop:shape-volume}, a contradiction.
\end{proof}

\begingroup
\small
\raggedright

\endgroup
\end{document}